\documentclass[12pt]{article}

\usepackage[T1]{fontenc}
\usepackage[utf8]{inputenc}
\usepackage[a4paper,margin=2.5cm]{geometry}
\usepackage{amsmath,amssymb}
\usepackage{graphicx}
\usepackage{float}
\usepackage{algorithm}
\usepackage{algpseudocode}
\usepackage{multicol}
\usepackage{xcolor}
\usepackage{soul}
\usepackage{textcomp}
\usepackage{url}
\usepackage{hyperref}
\hypersetup{hidelinks}
\usepackage{subfig}

\newcommand\numberthis{\addtocounter{equation}{1}\tag{\theequation}}
\renewcommand\hl{}

\title{Nonlinear Control Allocation:\\ A Polynomial Approach}

\author{Jan Belák\thanks{Faculty of Electrical Engineering, Czech Technical University in Prague, Department of Control Engineering. This person is the corresponding author of this paper and can be reached at belakjan@fel.cvut.cz.}, 
Martin Hromčík\thanks{Faculty of Electrical Engineering, Czech Technical University in Prague, Department of Control Engineering. He can be reached at xhromcik@fel.cvut.cz.} and 
Didier Henrion\thanks{LAAS-CNRS, CNRS, Université de Toulouse, France, and Faculty of Electrical Engineering, Czech Technical University in Prague, Department of Control Engineering. He can be reached at henrion@laas.fr.}}
\date{\today}

\begin{document}
\maketitle

\begin{abstract}
Control allocation   maps the virtual control inputs requested by a high-level control law onto the physical actuators of a system. When the actuator authority is strongly nonlinear, existing methods rely on local linearization, on incremental updates, or on general nonlinear programming, and may be inaccurate, slow, or trapped in local optima. This paper introduces Polynomial Control Allocation, in which the cost and constraints of the allocation problem are polynomial or rational functions of the states and inputs; radical, trigonometric and tabulated nonlinearities are accommodated by lifting. The first order optimality conditions then form a parametric system of polynomial equations. Its Gröbner basis and the associated multiplication matrices are computed once, offline. Online, the matrices are evaluated at the current state and all solutions are recovered from the eigenvectors of a random linear combination of them, so that the globally optimal allocation is selected at a predictable computational cost. The framework is illustrated on thrust control of a variable-pitch propeller and on multicopter force allocation, and its computational cost is quantified in terms of polynomial evaluation and eigenvector computation. The method is validated experimentally on a Quanser helicopter model at a 500 Hz control rate, where a two-layer allocation scheme combined with linear model predictive control decouples travel and elevation over multiple full rotations. Monte Carlo comparisons with nonlinear programming and incremental nonlinear control allocation show that Polynomial Control Allocation achieves the lowest allocation error together with the shortest computation time.
\end{abstract}

\section{Introduction}
Control allocation (CA) is a fundamental component of many aircraft \cite{matamoros2018incremental}, \cite{Burken2001}, \cite{10.2514/1.52234}, \cite{Sieberling2010} and spacecraft \cite{10.2514/1.G005834}, \cite{10.2514/1.G006266}, \cite{10.2514/1.G006915} control architectures. The CA system \hl{separates} the coupling interactions between states and inputs from higher-level control \hl{laws} \cite{le2022scalable}. \hl{For example, an aircraft airspeed controller requests propeller force, and the CA system calculates the propeller speed to achieve it. The influence of airspeed and air density on the force is resolved by the CA algorithm and is not included in the feedback controller design.} CA methods can also provide online reconfiguration in the case of actuator failure \cite{Burken2001}. The CA process applies to general nonlinear systems of the form
\begin{equation}
    \frac{d\vec{x}}{dt} = f(\vec{x})+ g(\vec{x},\vec{u}),
    \label{eq:gen_dyn_system}
\end{equation}
where $\vec{x}\in\mathbb{R}^n$ is the state vector, $\vec{u}\in\mathbb{R}^m$ is the control action vector, and the autonomous dynamics $f(\vec{x})$ are explicitly isolated. The nonlinear coupling of control action with system states, expressed in $g(\vec{x},\vec{u})$, complicates the design of control laws. This issue affects traditional nonlinear control approaches, such as feedback linearization and state inversion, which require the system to be affine with respect to control action. \hl{However, CA can be utilized to drive the system from Equation }\ref{eq:gen_dyn_system} \hl{to behave as the generic input-affine} \cite{matamoros2018incremental} form $\frac{d\vec{x}}{dt} = f(\vec{x})+ \hat{g}(\vec{x})\cdot \vec{v}$, \hl{from the perspective of higher-level control laws.} The $\vec{v} \in \mathbb{R}^n$ is a virtual control action vector, that is manipulated by the controller instead of $\vec{u}$.
%affecting the system state through the nonlinear function $\hat{g}(\vec{x})$.

In mechanical systems, the $\vec{u}$ variables represent real actuators (e.g. aircraft control surfaces), while $\vec{v}$ typically represents forces and moments for respective degrees of freedom as described in \cite{Durham1993} and \cite{10.2514/3.21201}. This system form enables the design of control laws using the virtual control action $\vec{v}$. The CA then assigns, or allocates, the control action $\vec{u}$ related to the prescribed virtual control $\vec{v}$ and system state $\vec{x}$ such that ideally the following set of algebraic equations
\begin{equation}
    \hat{g}(\vec{x})\cdot \vec{v} = g(\vec{x},\vec{u}),
    \label{eq:CA_definition}
\end{equation}
is solved. An expansion of this approach called Nonlinear Dynamic Inversion \cite{Sieberling2010} includes autonomous dynamics during allocation ($\hat{g}(\vec{x})\cdot \vec{v} = g(\vec{x},\vec{u})+f(\vec{x})$), which transforms the controlled dynamics into a set of first-order integrators. When the control action $\vec{u}$ is constrained, Equation \ref{eq:CA_definition} \hl{might not have a solution, and the CA method searches for an approximate solution using optimization techniques}. If the problem in Equation \ref{eq:CA_definition} constitutes a set of linear equations, we refer to linear CA. This problem is well understood and can be solved using generalized inverse methods in unconstrained cases \cite{Oppenheimer2006} and with linear/quadratic programming \hl{(LP/QP)} \cite{Burken2001}, \cite{Enns1998810} or daisy chaining \cite{Buffington1996} methods when control constraints are defined. In this paper, \hl{we do not consider linear CA techniques, instead }we consider $g(\vec{x},\vec{u})$ to be a set of multivariate polynomials. Such polynomial models are well established in the relevant domains of vehicle dynamics and control—consider, for example, NASA's Generic Transport Model and related flight stability analysis \cite{Chakraborty2011} or flight control validation and verification results \cite{Wagner2023}.

\section{Nonlinear Control Allocation Methods}
\label{sec:SOTA}
State-of-the-art nonlinear control allocation can be divided into three distinct approaches. The simplest approach, described in Section \ref{SOTA:loc_lin}, locally linearizes the function $g(\vec{x},\vec{u})$ for multiple values $\vec{x}_0$ and $\vec{u}_0$. Another approach, described in Section \ref{SOTA:incremental}, utilizes incremental control with online linearization of the control authority function. Situations with over- or under-actuated systems $(m\neq n)$ typically lead to optimization problems as described in Section \ref{SOTA:opt}.

\subsection{Piecewise Local Linearization Method}
\label{SOTA:loc_lin}
The traditional approach to CA locally linearizes the function $g(\vec{x},\vec{u})$ with respect to $\vec{u}$ at a nominal control vector $\vec{u}_0$. Substituting $\vec{x}$ with nominal state $\vec{x}_0$ yields
$\hat{g}(\vec{x}_0) = B = \frac{\partial}{\partial\vec{u}} g(\vec{x}_0,\vec{u}_0)$,
which is a constant-gain matrix. This limits the CA accuracy to a small region around the nominal values $\vec{x}_0$ and $\vec{u}_0$. When the control region exceeds the area accurately represented by linearization, the approach loses robustness and can lead to stability loss. \hl{An alternative is to linearize the function} $g(\vec{x},\vec{u})$ \hl{at each time step with current values} $\vec{x}_t$, $\vec{u}_t$  \hl{instead of} $\vec{x}_0$ \hl{and} $\vec{u}_0$ \hl{only, which moves the CA region in time, allowing gradual utilization of the entire control region. If the entire control region must be utilized in each control step, a nonlinear approach was developed in} \cite{Bolender2004}. It computes the matrix $B$ for multiple $\vec{u}_0$ values distributed across the control region, thereby approximating the nonlinear control authority function $g(\vec{x},\vec{u})$ using piecewise linear functions for given $\vec{x}_t$. This approach reformulates the control allocation problem as a piecewise linear program, which can be cast as a mixed-integer linear program, enabling exact solutions that better match the desired control action compared to traditional linear approaches. While \cite{Bolender2004} demonstrates the ability to solve nonlinear CA problems, it suffers from scalability issues, since mixed-integer linear programming is NP-hard.

\subsection{Incremental Nonlinear Control Allocation Method (INCA)}
\label{SOTA:incremental}
Another approach to solving the nonlinear CA problem utilizes online local linearization combined with incremental control \cite{Ostroff2002}. It employs a real-time updated Jacobian $B(\vec{x})$ of the control effectiveness function $g(\vec{x},\vec{u})$ to solve a linear control allocation problem despite nonlinear actuator authority. By handling incremental changes of control inputs instead of absolute values, the method can ensure accurate allocation as long as the maximum incremental step is sufficiently small. \hl{Another benefit is the straightforward implementation of rate limits, as they can be incorporated into saturation limits in each control step.} Consequently, the method appears most effective for densely sampled cases, %\hl{where controller frequency is sufficiently higher than control bandwidth of the system}, 
as discussed in Section \ref{sec:Performance comparison to State of the art}. See \cite{matamoros2018incremental} for an advanced flight control application. The main disadvantage of this method is its global non-optimality for non-convex problems. Similar to gradient-search methods, the Incremental Nonlinear Control Allocation method is sensitive to saddle points and local optima. Therefore, the method cannot easily escape such regions.

\subsection{Nonlinear Programming Control Allocation (NPCA)}
\label{SOTA:opt}
A generic optimization-based CA problem is described by the cost function $c(\vec{x},\vec{u})$, the constraint functions $l_1(\vec{x},\vec{u})$, $l_2(\vec{x},\vec{u})$, and the control variable limits $\vec{u}_{min}$, $\vec{u}_{max}$ as:
\begin{align*}
    \min &\;c(\vec{x},\vec{u})\\\numberthis \label{eq:generic_CA_problem}
    \text{s.t.}&\;l_1(\vec{x},\vec{u})=0&
    &\;l_2(\vec{x},\vec{u})\geq0&
    &\;\vec{u}_{min}\leq \vec{u}\leq \vec{u}_{max}
\end{align*}
Problem \ref{eq:generic_CA_problem}  can be reformulated as a system of algebraic equations (Karush-Kuhn-Tucker conditions \cite{Pedregal2004}), \hl{where} $\nabla_{\vec{u}}$ \hl{is gradient with respect to the vector} $\vec{u}$, $\lambda_i$ \hl{are equality- and} $\mu_j$ \hl{are inequality-constraint slack variables respectively:}
\begin{align*}
    &\nabla_{\vec{u}} c(\vec{x},\vec{u}) + \sum_i^m \lambda_i \cdot \nabla_{\vec{u}} l_{1,i}(\vec{x},\vec{u}) + \sum_j^k \mu_j \cdot \nabla_{\vec{u}} l_{2,j}(\vec{x},\vec{u}) =0& 
    \numberthis \label{eq:generic_CA_equations}
    &\forall i:  l_{1,i}(\vec{x},\vec{u}) = 0&
    &\forall j:  \mu_j\cdot l_{2,j}(\vec{x},\vec{u}) = 0,
\end{align*}
%5with control authority constraints expressed as
%\begin{align*}
%    (\vec{u}-\vec{u}_{max})&\leq0 &
%    (\vec{u}_{min}-\vec{u})&\leq0.
%\end{align*}
\begin{itemize}
    \item Computational complexity prevents usage for fast processes with high sampling rates
    \item Convergence is not guaranteed due to the gradient-descent nature of these methods
\end{itemize}

In this paper, we propose a novel polynomial CA method that enables the allocation of highly nonlinear and nonconvex problems without the limitations of the previously described methods.

\section{Polynomial Control Allocation (PCA)}
\label{sec:PCA}
Polynomial Control Allocation (PCA) utilizes the structure shown in Equation \ref{eq:generic_CA_problem}. The cost $c(\vec{x},\vec{u})$ and constraints $l_1$ and $l_2$ are expanded to rational functions $\frac{p(\vec{x},\vec{u})}{q(\vec{x},\vec{u})}$, where $p(\vec{x},\vec{u})$ and $q(\vec{x},\vec{u})$ are multivariate polynomial functions. Below, we indicate how certain non-polynomial objects can be included in the PCA framework.

\subsection*{Radical Functions}
If the problem includes a radical function, such as square root $\sqrt{x}$, the problem can be transformed to polynomial form by introducing an additional variable $\hat{x}$ and a constraint equation $\hat{x}^2 = x, \quad \hat{x} \geq 0$. For example, the system $F = \sqrt{x}\cdot u$ has the same solution set as the combined equations $F = \hat{x}\cdot u$ and $\hat{x}^2= x, \hat{x} \geq 0$.
\subsection*{Trigonometric Functions}
\label{subsec:trigonometric_substitution}
When a variable is represented only through trigonometric functions ($\sin(x)$, $\cos(x)$, or $\tan(x)$), the system can be precisely described using Euler's identity $\sin^2(x)+\cos^2(x) = 1$.
We introduce two additional variables for substitution $x_s = \sin(x)$ and $x_c = \cos(x)$ and add an equality constraint $x_s^2+x_c^2 -1 = 0$. \hl{This approach is demonstrated in the drone force control example from section} \ref{sec:Examples} \hl{and in the Outer CA Layer described in section }\ref{sec:ApplicationPCA}.
\subsection*{Rational equations}
\hl{Transformation of rational equations, represented by polynomials $f(\vec{x})$ and $g(\vec{x})$ and constant variable $a$ in form $ \frac{f(\vec{x})}{g(\vec{x})} = a$ can be transformed into polynomial equation $f(\vec{x}) = a\cdot g(\vec{x})$, where $g(\vec{x}) \neq 0$. We include $f(\vec{x}) = a\cdot g(\vec{x})$ form into the polynomial solver and check the results against the requirement $g(\vec{x}) \neq 0$ to select the correct solution.}
\subsection*{Other Analytical Functions}
When other non-polynomial functions are present in the problem, local approximation using Taylor series can be employed. Unlike traditional linearization procedures, the polynomial approach can use higher-order Taylor series to ensure robust results throughout the entire envelope.

\subsection*{Lookup Tables and Polynomial Interpolation}
One common approach to modeling highly nonlinear dynamics uses lookup tables, often populated with measured data points. An example would be aircraft aerodynamics captured from wind tunnel experiments. These tables can be approximated by low-order polynomials \hl{with respect to control variables, which are typically deflections of control surfaces. An example of aerodynamic authority of actuator} $\delta$ \hl{in roll axis} $ c_{l,\delta}$ \hl{is demonstrated in equation} \ref{eq:aircraft_actuators_func}.
\begin{equation}
    c_{l,\delta} = a(\alpha,\beta)+b(\alpha,\beta)\cdot\delta+c(\alpha,\beta)\cdot\delta^2+d(\alpha,\beta)\cdot\delta^3
\label{eq:aircraft_actuators_func}    
\end{equation}
Polynomial coefficients $a,b,c,d$ are defined as lookup tables that depend on the aircraft's aerodynamic angles $\alpha$ and $\beta$. The concept of polynomial modeling of aircraft dynamics was demonstrated by NASA's Generic Transport Model \cite{Chakraborty2011}. \hl{The approximation of aerodynamic data by low-order polynomials is utilized in the variable pitch propeller thrust control example from section} \ref{sec:Examples}.
\section{Gröbner-Basis Solver for Polynomial Control Allocation}
\label{sec:Grobner}
In this section, we introduce the polynomial control allocation (PCA) solver. We assume that the problem $F$ is a set of polynomial equations that describe the optimization \hl{equation} \ref{eq:generic_CA_equations}. To find solutions of $F$, we use Gröbner bases and numerical computations of the eigenvectors of multiplication matrices. A specific case of multiplication matrix for one-variable polynomials is the companion matrix often used for root evaluation\footnote{https://www.mathworks.com/help/matlab/ref/roots.html}. The Gröbner basis approach to polynomial optimization is classical and is explained in \cite{Coxc2005}, Chapter 2. It has already found applications in computer vision for fast estimation of relative or absolute camera pose from image correspondences \cite{Kukelova2012} and in parameterization of control systems \cite{Patil2014}. A more detailed description of the Gröbner basis can be found in Appendix \ref{appendix:GrobnerBasis}.

We propose constructing the Gröbner-basis PCA solver in two phases. Initially, a set of polynomial matrices is constructed offline. In the online phase, these matrices are evaluated for given state and input values and then subjected to eigenvector decomposition. This approach ensures real-time execution of the PCA solver.

\subsection{Offline Computations}
\subsection*{Gröbner-Basis Transformation}
The Gröbner basis is a structured set of polynomials in canonical form, allowing algorithmic solutions to polynomial algebraic equations through substitutions.
%For example, a set of polynomial equations
%\begin{align*}
%    x^2-y^2=(x-y)\cdot(x+y)&=0 & \numberthis \label{eq: Grobner-basis input example}
%    x^2+y&=0
%\end{align*}
%has a Gröbner-basis formulation
%\begin{align*}
% y^2 + y=y\cdot(y+1)&=0& \numberthis \label{eq: Grobner-basis example}
%    x^2 + y&=0.
%\end{align*}
%The roots of the first equation in \ref{eq: Grobner-basis example} are clearly $y=\{0,-1\}$, while the first equation in \ref{eq: Grobner-basis input example} has solutions $y=\{-x,x\}$ and requires further computation to obtain the solution.
Gröbner basis generation methods, such as Buchberger's Algorithm \cite{Conti1991}, rely on polynomial division and Gaussian elimination. These methods are computationally expensive and therefore ill-suited for real-time applications. However, we propose precomputing the Gröbner basis of the PCA problem offline, thus avoiding real-time performance concerns.

\subsection*{Generation of Multiplication Matrices}
An effective way to find the solution of a Gröbner-basis set of polynomial equations F is through multiplication matrices. A multiplication matrix is a finite-dimensional matrix that represents how multiplying by a variable transforms Gröbner-basis monomials related to its system of equations. The matrix eigenvalues correspond to the roots of the Gröbner-basis system, which forms the basis for the online computations of our method. A multiplication matrix for a Gröbner-basis system is computed as follows:

\begin{algorithm}
\caption{Multiplication matrix generation for \hl{the set of polynomial equations F}}\label{alg:multiplication matrix}
\begin{algorithmic}[1]
\Require Equations F are polynomial w.r.t variables.
\Require MAPLE software (not implemented in MATLAB)
\State Select polynomial order (MAPLE "tdeg")
\State Generate Gröbner-basis G (MAPLE "Basis")
\State Find monomial basis for G (MAPLE "NormalSet")
\ForAll{optimization variables $x_i$}
    \State Compute multiplication matrix $M_{x_i}$ (MAPLE "MultiplicationMatrix")
    \State Export $M_{x_i}$ to csv file (MAPLE "ExportMatrix")
    \State Import $M_{x_i}$ into MATLAB as sparse symbolic matrix
\EndFor
\end{algorithmic}
\end{algorithm}
%Note that the size of the multiplication matrix corresponds to the number of solutions to the system.
\subsubsection*{Simple 1 variable example}
\hl{The generation of a multiplication matrix for %an ideal I (see definition in appendix} \ref{appendix:GrobnerBasis})\hl{, generated by 
a polynomial system of equations $F$, assuming $c\neq0$. The Groebner basis $G$ is identical to  $F$ as $F$ is a single polynomial.}
\begin{align*}
    F = G = a + b\cdot{x}+ c\cdot{x}^2 &=0 &c\neq&0
\end{align*}
\hl{First, we define the monomial basis of the system as $m = \{1,x\}$.%
%quotient algebra $R[x] / I$: $m = \{1,x\}$.
We find the multiplication matrix of $x$ by multiplying the monomial basis by the variable $x$ and expressing the results using the Groebner basis G (which is identical to F in this case):}
\begin{itemize}
    \item $1\cdot x = x$: $x$ is part of the monomial basis and does not have to be computed from $F$
    \item $x\cdot x = x^2 = \frac{-a-b\cdot x}{c}$: $x^2$ is not part of $m$ and must be expressed from equation $F$ as $c\cdot x^2 = -a-b\cdot x$
\end{itemize}
\hl{We express the results using $m$ vector coordinates as}
\begin{align*}
    x &= \begin{bmatrix}
         0 & 1
     \end{bmatrix}\cdot \begin{bmatrix}
         1\\x
     \end{bmatrix} & \frac{-a-b\cdot x}{c}&= \begin{bmatrix}
         -\frac{a}{c} -\frac{b}{c}
     \end{bmatrix}\cdot \begin{bmatrix}
         1\\x
     \end{bmatrix}.
     \end{align*}
    \hl{The multiplication matrix is then created by stacking the two vectors together into 2$\times$2 matrix} 
     \begin{equation*}
         M_x = \begin{bmatrix}
          0 & 1\\
          -\frac{a}{c} &-\frac{b}{c}.
     \end{bmatrix}
     \end{equation*}
     \hl{The eigenvalues of this matrix correspond to the solutions of the original polynomial. This example shows only the basic ideas behind generating a multiplication matrix. A more detailed description can be found in }\cite{MutliplicationMatrixExpl}.
\subsection{Online Computations}
The steps needed for real-time computation of the PCA solution are shown in Algorithm \ref{alg:PCA_solver}, with individual steps explained in the following paragraphs.

\subsection*{Evaluation of Multiplication Matrices}
The first step of real-time computation is evaluating rational multiplication matrices for current values of system states and inputs. Evaluation of numerous polynomials can be computationally expensive and can significantly limit achievable sampling rates.

\subsection*{Computation of the System's Eigenvectors}
With all multiplication matrices evaluated, a generalized matrix is constructed as a randomized sum of multiplication matrices $M = \sum_{i=1}^k a_i\cdot M_{x_i}$, where $M_{x_i}$ is the multiplication matrix of a given variable and $a_i\in(0,1)$ is a random coefficient, with $\sum_{i=1}^k a_i= 1$ and $k$ being the number of multiplication matrices. $k$ also corresponds to the number of variables in the Gröbner basis.

\subsection*{Computation of the Solution}
To compute the solution, we find the eigenvectors $V$ of matrix $M$. The solution to a specific variable $x_i$ is then 
\begin{equation}
    x_i = \text{diag}(V^T\cdot M_{x_i}\cdot V).
\end{equation}
As can be seen, the complexity of online computation is dictated by the evaluation of multiplication matrices and the computation of eigenvectors. Computational complexity is discussed in detail in Section \ref{sec:Performance}.

\subsection{Generalized Eigenvalue Algorithm}
To calculate eigenvectors and eigenvalues, we utilize the generalized Schur decomposition. %This method is implemented by default in MATLAB's `eig` function and performs well for generic (non-symmetric) matrices. 
This algorithm has been shown to be backward stable for all practical purposes \cite{Chandrasekaran1996} for a large class of structured matrices. In our experience, the resulting accuracy of eigenvalue calculations is at the level of system precision for all considered PCA cases.% However, we recommend properly scaling the variables of the PCA problem, meaning that all parameter ranges are similar. %We also recommend verifying the calculated results through the original problem to ensure the accuracy of the solution.

\begin{algorithm}
\caption{Online computation of PCA solution}\label{alg:PCA_solver}
\begin{algorithmic}[1]
\ForAll{Multiplication matrices $M_{x_i}$}
    \State Evaluate the symbolic $M_{x_i}$ matrix for current non-optimized variables
\EndFor
\State Calculate $M= \sum_i \alpha_i M_{x_i}$, where $\sum_i \alpha_i = 1$
\State Calculate eigenvectors $V$ of matrix $M$
\ForAll{Multiplication matrices $M_{x_i}$}
    \State Calculate vector of solutions $\vec{x}_i=  \text{diag}( V^T \cdot M_{x_i} \cdot V)$
\EndFor
\State \hl{Select only real solution (with zero imaginary value: }$\forall i: \Im(x_{i})=0$), compliant with KKT conditions
\State Verify results using the input model
\end{algorithmic}
\end{algorithm}

\section{PCA Case Studies}
\label{sec:Examples}
In this section, we present selected practical problems that can be effectively addressed by our PCA solution. Essential ideas are elaborated here, while numerical solutions can be found in the dedicated git repository at \cite{iuITwfwC9BcxK0Lg}.

\subsection*{Variable Pitch Propeller Thrust Control: Multidimensional High-Order Actuation Function}
\begin{figure}[h]
    \centering
    \subfloat[\centering $c_T(J)$]{{\includegraphics[width=.54\linewidth]{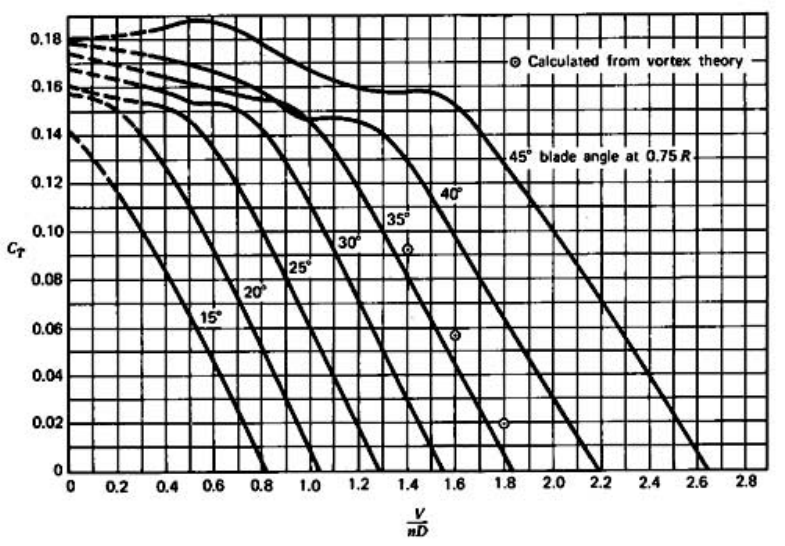} }}%
    \qquad
    \subfloat[\centering $c_P(J)$]{{\includegraphics[width=.40\linewidth]{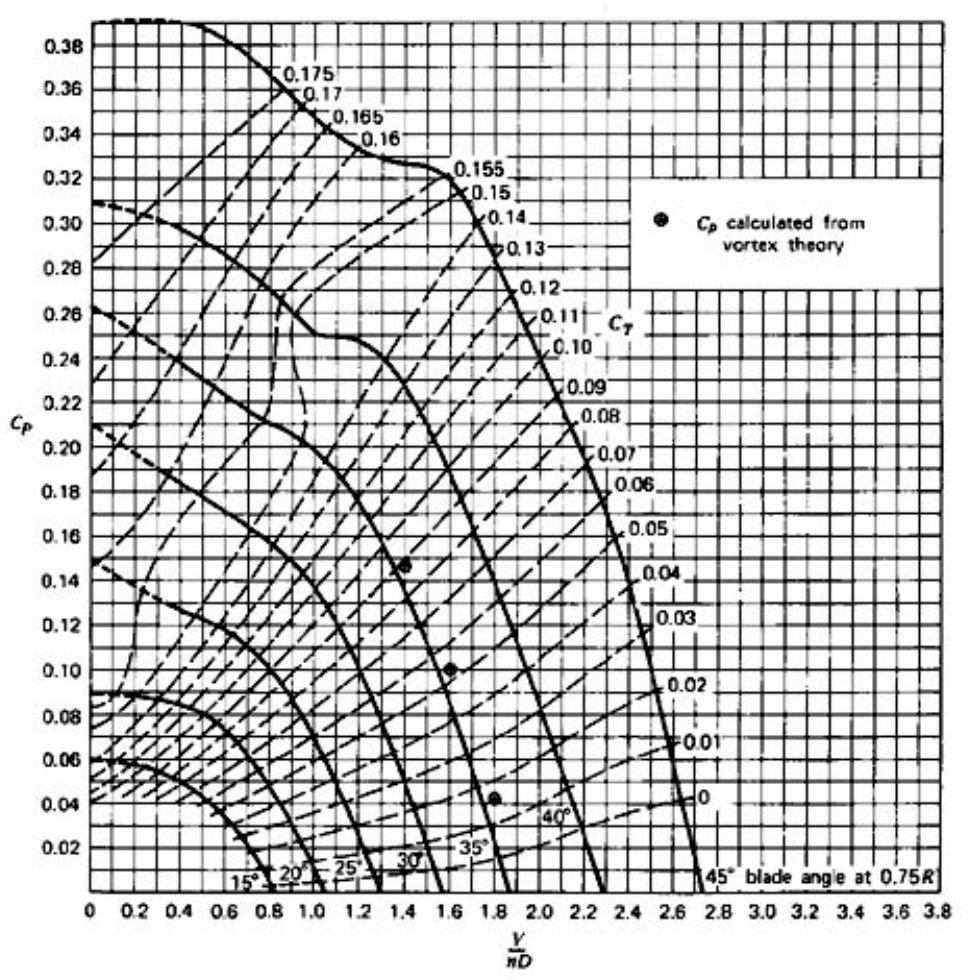} }}%
    \caption{Propeller thrust and power coefficients, reproduced from \cite{McCormick1979}}
    \label{fig:k_T}
\end{figure}

Modern electric engines, utilized in drones and potentially in electrically powered aircraft, can operate across a wide range of rotational speeds with nearly constant efficiency. This motivates efforts to simultaneously control the propeller pitch angle $\beta$ and its rotational speed $n$. With two independent variables, we can optimize control to maximize efficiency. We describe the propeller thrust $T$ as $T = n^2\cdot c_T\cdot\rho$, \hl{where $\rho$ is air density} and $c_T$ is the thrust coefficient, dependent on the pitch angle, vehicle speed, and propeller rotational speed, depicted in %\cite{McCormick1979}. The relation is depicted in 
Figure \ref{fig:k_T}.
For $c_T < 0.14$, we can approximate the relation as a linear map of $\beta$ and the advance ratio $J = v/(n\cdot D)$,\hl{ dependent on airspeed $v$ and propeller diameter $D$,} as $c_T \approx a+  b\cdot J + c\cdot \beta$, where $a$, $b$, $c$ are suitably selected parameters. \hl{Combining these relations gives us the final form of the thrust equation:}
\begin{equation}
   T(n,\beta,v,\rho) =   n^2\cdot (a+ b\cdot\frac{v}{n\cdot D} + c\cdot \beta)\cdot\rho 
\end{equation}

The propeller power $P =  \rho \cdot c_P\cdot n^3$ is expressed through the power coefficient $c_P$ and propeller speed $n$ \cite{McCormick1979}. The coefficient $c_P$ is nonlinearly dependent on $\beta$ and $J$ as %depicted in \cite{McCormick1979} 
shown in Figure \ref{fig:k_T}, and we approximate it as $c_P \approx J^2\cdot p_{20} + \beta \cdot p_{01} + p_{00}$, where $p_{ij}$ are suitably selected parameters. \hl{The expanded form of the equation is}
\begin{equation}
    P(n,\beta,v,\rho) = \rho \cdot (\frac{v^2}{n^2\cdot D^2}\cdot p_{20} + \beta \cdot p_{01} + p_{00})\cdot n^3
\end{equation}

\hl{The system must track the required thrust $T_{ref}$, set by the feedback controller, and comply with structural limits on $\beta$ and $n$. From these relations, we can build the optimization problem in the form described by equation }\ref{eq:generic_CA_problem}\hl{, minimizing power $P$ for given thrust $T_{ref}$ by optimizing propeller speed $n$ (replaceable by the advance ratio $J$ as they are inversely related) and pitch angle $\beta$}
\begin{align*}
    \min_{n,\beta} &\;P(n,\beta,v,\rho)\\
    \text{s.t.}&\; T_{\text{ref}} =T(n,\beta,v,\rho)&
    &\; \beta_{\min}-\beta \leq 0 &\beta-\beta_{\max}\leq0&
    &\; n_{\min} -n\leq& 0 & n- n_{\max}\leq0.
\end{align*}
\hl{Comparing this optimization task to the general definition from equation} \ref{eq:generic_CA_problem}, \hl{we see that $c(\vec{x},\vec{u}) = P(n,\beta,v,\rho)$, $l_1 = T_{\text{ref}}-T(n,\beta,v,\rho)$, and $l_2 = \begin{bmatrix}
    &\; \beta_{\min}-\beta &\beta-\beta_{\max} & n_{\min} -n & n- n_{\max}
\end{bmatrix}$. The state variables $\vec{x}$ are $\rho$ and $v$ while control variables $\vec{u}$ are $n$ and $\beta$.} Transformation of the optimization task into a set of equations according to equation \ref{eq:generic_CA_equations} results in \hl{the following polynomial system F with slack variables $\mu_1$ to $\mu_4$ represent inequality constraints for $\beta$ and $n$ limits, while $\lambda_T$ variable represents the equality constraint of thrust}.
\begin{align*}
    \frac{2 \rho  p_{20} \,v^{2}}{D^{2}}-3 \rho  (p_{00} +\frac{p_{20} \,v^{2}}{n^{2} D^{2}}+p_{01} \beta ) n^{2}-\lambda_T (2 D c n \beta \rho +2 D a n \rho +b \rho  v )-\mu_{1} +\mu_{2}&=0\\
    \end{align*}
    \vspace{-70pt}
    \begin{align*}
    -D c \lambda_T \,n^{2} \rho -n^{3} p_{01} \rho -\mu_{3} +\mu_{4}&=0&
    -D c \,n^{2} \beta \rho -D a \,n^{2} \rho -b n \rho  v +D T_{ref}&=0\\
    \end{align*}
    \vspace{-70pt}
    \begin{align*}
    \mu_1 (-n + n_{\min}) &=0 &
    \mu_2 (-n_{\max} +n)&=0 &
    \mu_3 (-\beta +\beta_{\min} )&=0 &
    \mu_4 (-\beta_{\max} +\beta )&=0
\end{align*}
The solution to the system of equations is in the form of multiplication matrices and is provided in the related git repository \cite{iuITwfwC9BcxK0Lg}, as the individual relations are too large to display here. The seven resulting multiplication matrices, corresponding to two optimization variables and five constraint-related variables, are $7\times7$. Selecting the correct solution requires constraint verification through selection of positive semi-definite solutions for all constraint-related variables.

\subsection*{Drone Force Control: Nonlinear Multidimensional Optimization}
Maneuvering of multi-copters in the environment can be achieved through force allocation via thrust and attitude configuration, as discussed in \cite{Stephan2025}. The force equations in the global coordinate frame are
\begin{align*}
    F_{xe} &= T\cdot \sin(\theta) \cdot\cos(\psi)&
    F_{ye} &= T\cdot \cos(\theta) \cdot\sin(\psi)&
    F_{ze} &= T\cdot \cos(\theta) \cdot\cos(\psi),
\end{align*}
where $T$ is drone thrust, while $\theta$ and $\psi$ are the drone roll and pitch angles, respectively. We substitute the trigonometric functions with the equality constraint from Section \ref{subsec:trigonometric_substitution}. The resulting form is
\begin{align*}
    F_{xe} &= T\cdot \theta_s \cdot\psi_c &
    F_{ye} &= T\cdot \theta_c \cdot\psi_s &
    F_{ze} &= T\cdot \theta_c \cdot\psi_c& 
    1&=\theta_c^2+\theta_s^2 &
    1&=\psi_c^2+\psi_s^2
\end{align*}

During low-altitude maneuvering, the designer can require precise vertical force control by limiting $F_{xe}$ and $F_{ye}$ when necessary. \hl{This approach ensures precise tracking o $F_{ze}$ by reducing tracking accuracy of $F_{xe}$ and $F_{ye}$. To achieve the desired direction of the resulting force in the $x$-$y$ plane, we add a requirement that their relative error $\frac{F_{xe}}{F_{xe,\text{ref}}} = \frac{F_{ye}}{F_{ye,\text{ref}}} = 1-\Lambda$ is identical. Minimizing $\Lambda^2$ ensures the best possible tracking of $F_x$ and $F_y$ forces while enforcing precise $F_z$ tracking.} We describe the corresponding CA problem according to Equation \ref{eq:generic_CA_problem} as
\begin{align*}
    \min &\; \Lambda^2\\
    \text{s.t.}&\; F_{ze}-F_{ze,\text{ref}}=0&
    &\; F_{xe}-(1-\Lambda)\cdot F_{xe,\text{ref}}=0&
    &\; F_{ye}-(1-\Lambda)\cdot F_{ye,\text{ref}}=0\\
    &\;\theta_c^2+\theta_s^2-1=0&
    &\;\psi_c^2+\psi_s^2-1=0&
    &\;T-T_{\max}\leq0
\end{align*}

The problem was transformed to an equivalent equality problem, defined in the system of equations \ref{eq:generic_CA_equations}. The optimization variables are $\Lambda$, $T$, $\theta_c$, $\theta_s$, $\psi_c$, $\psi_s$, and six constraint-related variables. The corresponding Gröbner basis is too large to display in this work, with 11 multiplication matrices. The multiplication matrices are $11\times11$, yielding 11 complex solutions to the system of equations. From the solutions, we select only real solutions with constraint variables being positive semi-definite. Details of the solution is provided in the accompanying MAPLE script.
\section{Computational Complexity Evaluation}
\label{sec:Performance}
We focus on the online portion of the algorithm. Computational complexity is dictated by two algorithms:
\begin{itemize}
    \item Evaluation of polynomials in multiplication matrices
    \item Computation of eigenvectors of the generalized matrix $M$
\end{itemize}

The remaining computations, such as matrix multiplications, are several orders of magnitude faster and are not considered in this analysis. To estimate computational requirements, we utilize Monte Carlo simulations with random inputs on a typical personal computer \footnote{Lenovo ThinkPad with Intel Core I7-10750H CPU and 16 GB of RAM memory, NVIDIA Quadro T1000 Max-Q 4GB, Windows 10, Matlab 2024b}. We assume that the overall computational requirement is the sum of individual tasks computed in sequence. Each task is evaluated individually to investigate its influence on the overall complexity.

\subsection*{Evaluation of Polynomial Coefficients in Multiplication Matrices}
Each coefficient of the multiplication matrix is a polynomial that must be evaluated in each control step. The matrices are often sparse; however, the order of the multiplication matrix polynomials significantly exceeds the order of the original polynomial CA problem from equation \ref{eq:generic_CA_problem}.

\begin{figure}[H]
    \centering
    \includegraphics[width=0.7\linewidth]{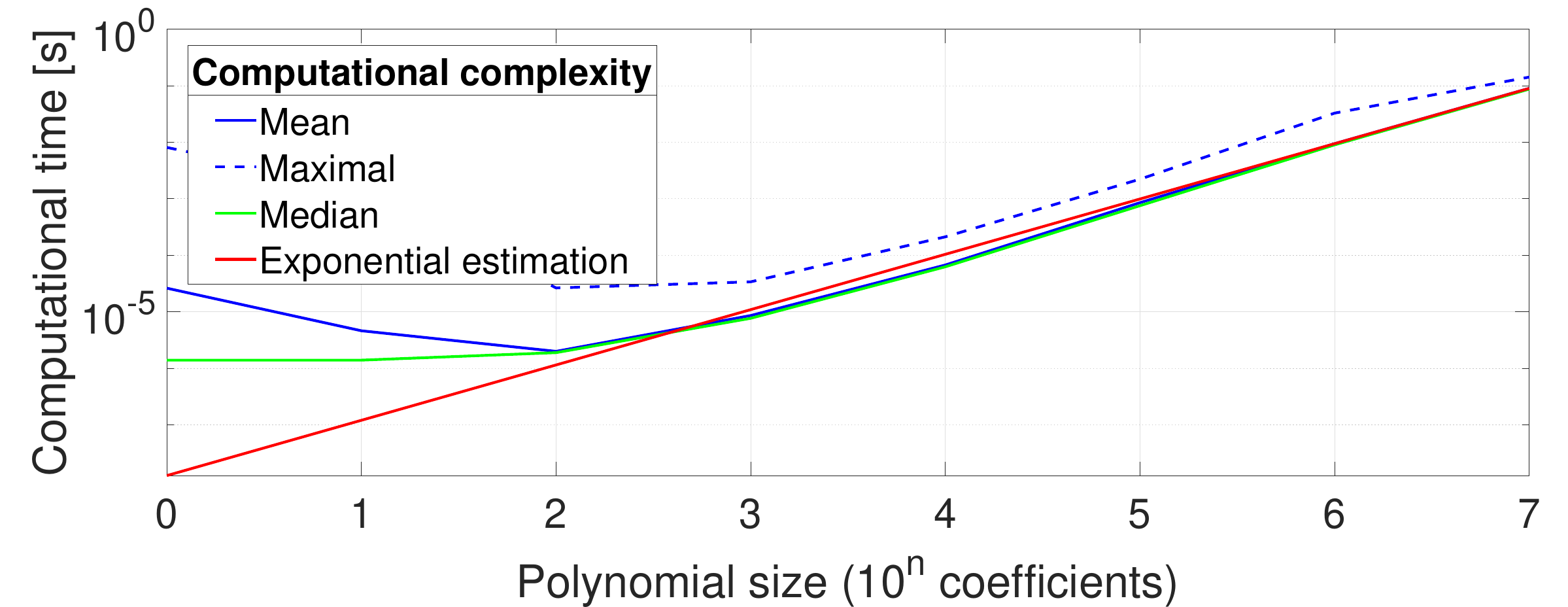}
    \caption{Time requirements for evaluating a polynomial with $10^n$ coefficients}
    \label{fig:time_poly}
\end{figure}

To provide a PCA problem designer with a practical guideline for the resulting computational complexity, we determined the computation time required to evaluate polynomials of varying orders. The relation between the number of coefficients of the multiplication matrix (corresponding to polynomial order and matrix size) and the time needed to evaluate the matrix is shown in Figure \ref{fig:time_poly}. We see that the system can calculate up to $10^5$ coefficients in about 1 ms. Due to exponential growth in time complexity, in 100 ms it can compute only about $10^3$ coefficients.

The number of coefficients of the Section \ref{sec:Examples} use cases and the corresponding computational times $\tau_{poly}$ are:
\begin{itemize}
\begin{multicols}{3}
    \item Fuel cell $\approx 10^1$: $\tau_{poly}\approx 10^{-2}$ ms
    \item Drone $\approx 10^4$: $\tau_{poly} \approx 10^{-1}$ ms
    \item Propeller $\approx 10^5$: $\tau_{poly} \approx 1$ ms
\end{multicols}
\end{itemize}
\subsection*{Computation of Eigenvectors of the Generalized Matrix}
\begin{figure}[H]
    \centering
    \includegraphics[width=0.7\linewidth]{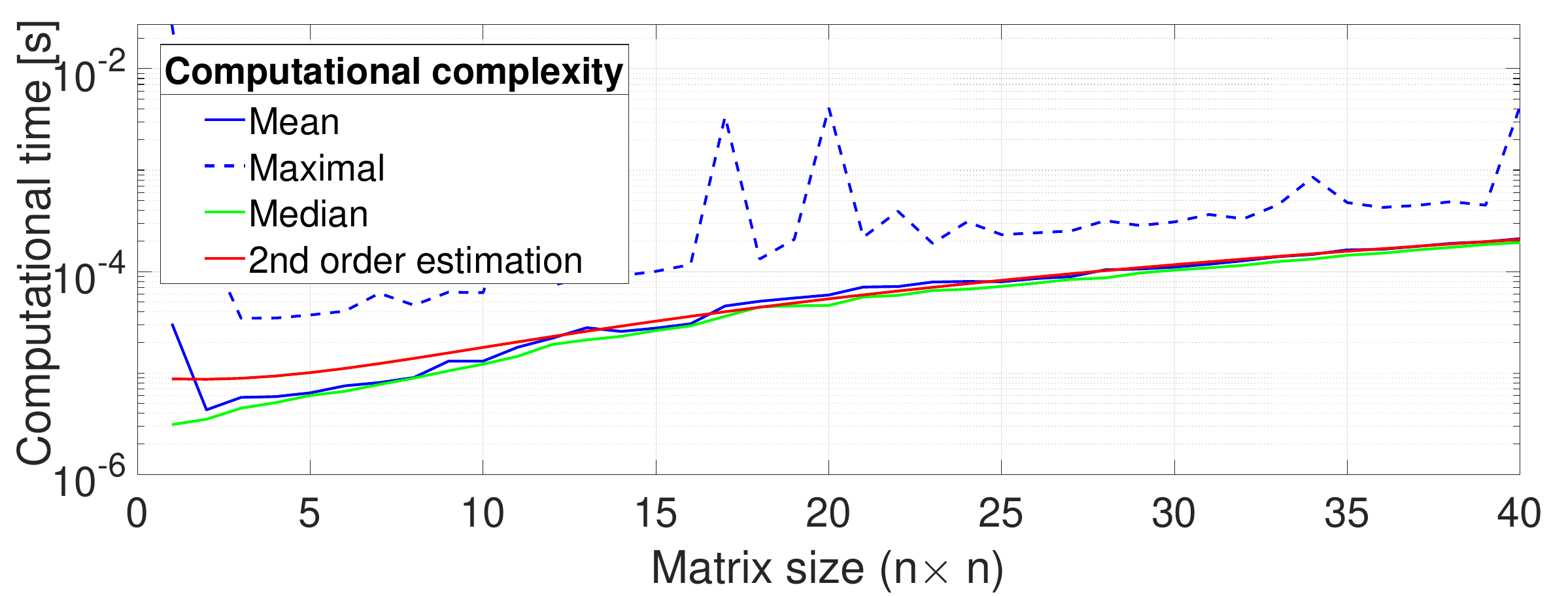}
    \caption{Time requirements for computing eigenvectors of matrices with varying size}
    \label{fig:time_eig}
\end{figure}

The method is based on computing eigenvectors of the linear combination of all multiplication matrices. The size of this general matrix corresponds to the number of solutions the problem has. We can estimate computational time for different matrix sizes with random coefficients. For each size, we tested $10^3$ different random matrices to estimate the trend between matrix size and computation time. Results are shown in Figure \ref{fig:time_eig}. As can be seen, computational complexity corresponds to $O(n^2)$, where $n$ is the matrix size, which aligns with the reported estimation $O(n^2)$ for structured matrices \cite{srOhbTQCDKdrMKYk}. Therefore, in most scenarios, we can expect that real-time applicability will not be limited by eigenvector computation. Specifically, the matrix size of the use cases and their computational times $\tau_{eig}$ are:
\begin{itemize}
\begin{multicols}{3}
    \item Fuel cell  $n=5$: $\tau_{eig} \approx 10^{-2}$ ms
    \item Drone $n=12$: $\tau_{eig} \approx 10^{-2}$ ms
    \item Propeller $n=7$: $\tau_{eig} \approx 10^{-2}$ ms
\end{multicols}
\end{itemize}
\hl{The overall achievable sampling frequency is determined by the sum of corresponding $\tau_{eig}$ and $\tau_{poly}$. We observe that more complex problems will likely be limited by the polynomial evaluation and not by the eigenvector computation.}
\section{Experimental Verification of PCA: Quanser 3DOF Helicopter Control}
\label{sec:Quanser}
\begin{figure}[h]
    \centering
    \includegraphics[width=.5\linewidth]{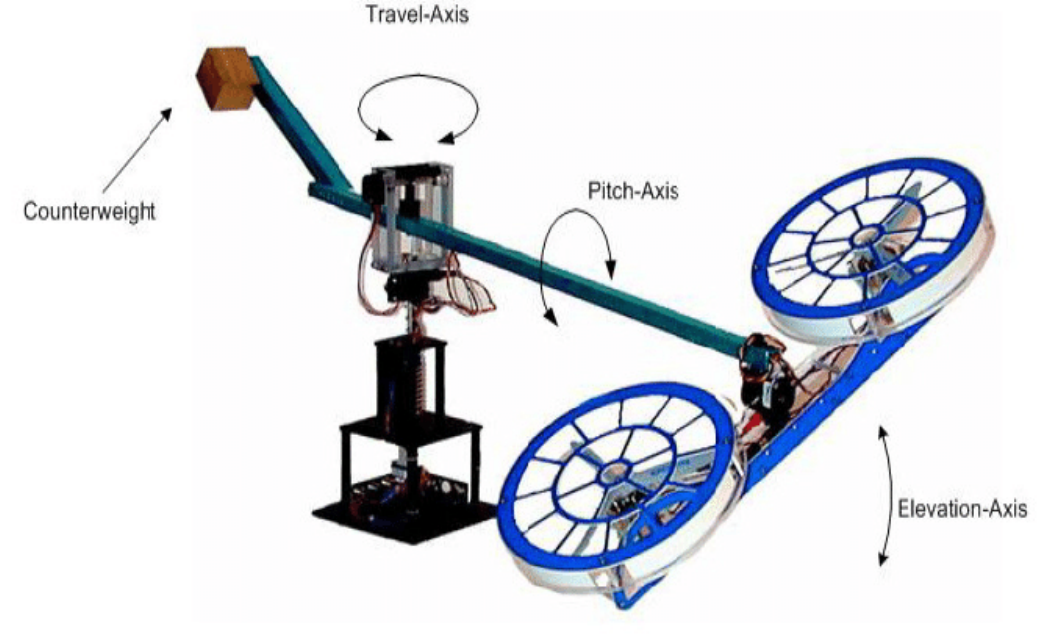}
    \caption{Quanser 3DOF Helicopter system}
    \label{fig:Quanser}
\end{figure}

To verify the newly developed polynomial control allocation method, we performed real-world experiments on a widely available and well-understood Quanser 3DOF helicopter system\footnote{Available at https://www.quanser.com/products/3-dof-helicopter/}, shown in Figure \ref{fig:Quanser}. The system was selected for its availability and nonlinear control action authority.

\subsection{Control Interface}
The Quanser system is connected to a standard computer with communication interfaces.% showcased in Figure \ref{fig:Quanser_interface}. 
The system measures all three angular states—pitch angle $\theta$, elevation angle $\epsilon$, and travel angle $\lambda$—using incremental sensors. The controller commands voltages for two propeller engines. The control frequency is 500 Hz, which is considered representative of challenging aerospace control.

%\begin{figure}[h]
%    \centering
%    \includegraphics[width=1\linewidth]{QuanserInterface.png}
%    \caption{Control interface of Quanser system}
%    \label{fig:Quanser_interface}
%\end{figure}

\subsection{Mathematical Model}
An experimentally verified nonlinear model for the Quanser system was developed in \cite{kf4EALo4Jr2YU6EW}, \hl{equation (5e)}. In this paper, we list only the dynamical equations and modifications performed to enable the polynomial CA. The system angular dynamics are
\begin{align*}
    \ddot{\lambda} &= -\text{sign}(\dot{\lambda})\cdot c_\lambda\cdot |\dot{\lambda}|^{2/3}- b_\lambda\cdot \cos(\epsilon)\cdot\sin(\theta)\cdot f_s\\
    \ddot{\theta} &= -a_\theta \cdot \cos(\epsilon) \cdot \sin(\theta)-c_\theta\cdot \dot{\theta} + b_\theta \cdot f_d\\
    \ddot{\epsilon} &= a_{\epsilon_1}\cdot\sin(\epsilon-\epsilon_0) -a_{\epsilon_2}\cdot \sin(\epsilon)\cdot \cos(\theta)-\text{sign}(\dot{\epsilon})\cdot c_{\epsilon}\cdot |\dot{\epsilon}|^{2/3}+b_\epsilon \cdot\cos(\theta) \cdot f_s,\\
\end{align*}
where $a_\theta$, $a_{\epsilon_1}$, $a_{\epsilon_2}$, $b_{\lambda}$, $b_{\epsilon}$, $b_{\theta}$, $c_\lambda$, $c_\epsilon$, and $c_\theta$ are system coefficients. The forces generated by propellers are included through relations $f_s = F_r(u_r) + F_l(u_l)$ and $f_d = F_r(u_r) - F_l(u_l)$, where $F_r(u_r)$ and $F_l(u_l)$ are forces generated by the right and left propeller, respectively. The forces $F$ depend on voltages $u$ through a third-order polynomial $F = p_1 u + p_3 u^3$, where $p_1$ and $p_3$ are polynomial coefficients. The third-order polynomial differs from the second-order polynomial proposed in the original paper \cite{kf4EALo4Jr2YU6EW}. These modifications enable to cover negative thrust values. %However, the authors conducted experiments mainly with positive thrust. Consequently, coefficients $p_1$ and $p_3$ were fitted to the positive thrust area and do not reflect the different function shape for negative voltages. The model can handle negative voltages (and thrust) but with reduced accuracy.

\subsection{Application of Polynomial Control Allocation Method}
\label{sec:ApplicationPCA}
The model can be rewritten in the form of Equation \ref{eq:gen_dyn_system} through definitions of autonomous dynamics and control authority function, assuming $\vec{x} = [\dot{\lambda},\dot{\epsilon},\dot{\theta}]^T$ and $\vec{u} = [u_r,u_l]^T$, in \hl{the} form $\frac{d\vec{x}}{dt} = f(\vec{x})+g(\vec{x},\vec{u})$ as
\begin{equation*} 
     \underbrace{\begin{bmatrix}
        -\text{sign}(\dot{\lambda})\cdot c_\lambda\cdot |\dot{\lambda}|^{2/3}\\
        a_{\epsilon_1}\cdot\sin(\epsilon-\epsilon_0) -a_{\epsilon_2}\cdot \sin(\epsilon)\cdot \cos(\theta)-\text{sign}(\dot{\epsilon})\cdot c_{\epsilon}\cdot |\dot{\epsilon}|^{2/3}\\
        -a_\theta \cdot \cos(\epsilon) \cdot \sin(\theta)-c_\theta\cdot \dot{\theta}
    \end{bmatrix}}_{f(\vec{x})}+ \underbrace{\begin{bmatrix}
        -b_\lambda\cdot \cos(\epsilon)\cdot\sin(\theta)\cdot (p_1\cdot u_r +p_3\cdot u_r^3+p_1\cdot u_l +p_3\cdot u_l^3)\\
        b_\epsilon \cos(\theta)\cdot (p_1\cdot u_r +p_3\cdot u_r^3+p_1\cdot u_l +p_3\cdot u_l^3)\\
        b_\theta \cdot (p_1\cdot u_r +p_3\cdot u_r^3-p_1\cdot u_l -p_3\cdot u_l^3)
    \end{bmatrix}}_{g(\vec{x},\vec{u})}
\end{equation*}
For the control system to independently address the travel angle $\lambda$ and the elevation angle $\epsilon$, the system must be able to control and track the pitch angle $\theta$. The angle $\theta$ is therefore considered dependent and will be controlled by the inner CA layer of the control system, combined with the allocation of voltages according to the overall desired thrust. The elevation and travel angles depend on the sum of thrust $f_s$ and pitch angle $\theta$, which will be the control variables for the outer CA layer.

\subsubsection*{Inner CA Layer}
The inner layer allocates propeller voltages to track the required pitch acceleration $\ddot{\theta}$ and overall thrust vector $f_s$. The thrust command is considered most crucial, as it is linked to elevation tracking\hl{ and must be tracked precisely. The pitch acceleration $\ddot{\theta}$ is secondary, so its tracking error $\Lambda_\theta$ is allowed but must be minimized.} The optimization problem is written as \hl{minimizing pitch acceleration error $\Lambda_\theta^2$ while tracking $f_{s,ref}$ in the constrained region of $u_l$ and $u_r$}
\begin{align*}
    \min_{u_l, u_r}\;& \Lambda_\theta^2 \numberthis\label{eq:InnerCA}\\
    \text{s.t.}\;& \ddot{\theta}= -a_\theta \cdot \cos(\epsilon) \cdot \sin(\theta)-c_\theta\cdot \dot{\theta}+b_\theta \cdot \underbrace{(p_1\cdot u_r +p_3\cdot u_r^3-p_1\cdot u_l -p_3\cdot u_l^3)}_{f_d}=(1-\Lambda_\theta)\cdot \ddot{\theta}_{\text{ref}}\\
    &f_s =p_1\cdot u_r +p_3\cdot u_r^3+p_1\cdot u_l +p_3\cdot u_l^3=f_{s,\text{ref}}\\
    &u_{\min}\leq u_l\leq u_{\max} \quad u_{\min}\leq u_r\leq u_{\max}.
\end{align*}

The problem is polynomial with respect to optimization variables $u_l$ and $u_r$. Consequently, we can solve this optimization problem with the proposed polynomial control allocation method. The generation of first-order derivatives and multiplication matrices is documented in our git repository \cite{iuITwfwC9BcxK0Lg}. The resulting multiplication matrices have a size of $n=33$, which should still allow for computation time below 1 ms on the considered hardware (refer to Section \ref{sec:Performance}).

\subsubsection*{Outer CA Layer}
The outer layer allocates the desired pitch angle $\theta$ and overall thrust $f_s$ to track the required travel acceleration $\ddot{\lambda}$ and elevation acceleration $\ddot{\epsilon}$. \hl{Precise elevation tracking} is desired to reduce elevation divergence during fast changes in travel angle.\hl{ The travel acceleration error $\Lambda_\lambda$ is allowed, but to be minimized.} We write the optimization problem as follows:
\begin{align*}
    \min_{f_s, \theta}\;& \Lambda_\lambda^2 \numberthis\label{eq:OuterCA}\\
    \text{s.t.}\;& \ddot{\lambda}=-\text{sign}(\dot{\lambda})\cdot c_\lambda\cdot |\dot{\lambda}|^{2/3}-b_\lambda\cdot \cos(\epsilon)\cdot\sin(\theta)\cdot f_s=(1-\Lambda_\lambda)\cdot \ddot{\lambda}_{\text{ref}}\\
    &\ddot{\epsilon}= a_{\epsilon_1}\cdot\sin(\epsilon-\epsilon_0) -a_{\epsilon_2}\cdot \sin(\epsilon)\cdot \cos(\theta)-\text{sign}(\dot{\epsilon})\cdot c_{\epsilon}\cdot |\dot{\epsilon}|^{2/3}+ b_\epsilon \cos(\theta)\cdot f_s =\ddot{\epsilon}_{\text{ref}}\\
    &f_{s,\min}\leq f_s\leq f_{s,\max}
\end{align*}

The problem with this formulation lies in the non-polynomial expression of $\theta$. Accelerations $\ddot{\lambda}$ and $\ddot{\epsilon}$ depend on the sine and cosine of angle $\theta$. We reformulate the problem by creating two additional variables and one constraint relation:
\begin{align}
c_\theta &= \cos (\theta) & s_\theta &= \sin (\theta) & c_\theta^2+s_\theta^2 &= 1.
\end{align}

The modified optimization problem is polynomial with respect to optimization variables and can utilize the developed PCA method. As with the inner layer, we do not include the generation of multiplication matrices in this paper and provide a link to git with the generated matrices. The matrices have size $n = 9$, which is much smaller than the inner layer and is therefore not assumed to be computationally critical.

\subsection{Feedback Control Law}
\begin{figure}[H]
    \centering
    \includegraphics[width=0.7\linewidth]{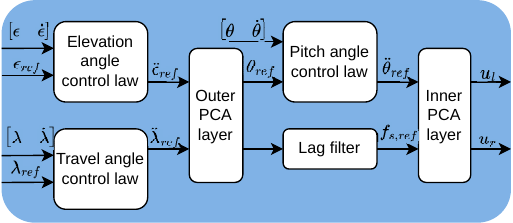}
    \caption{Control system of the Quanser system}
    \label{fig:ControlScheme}
\end{figure}

The system is controlled through a typical position-attitude hierarchical control scheme (see \cite{Nguyen2020}) with external tracking commands $\epsilon$ and $\lambda$ and internal tracking commands $\theta$, generated by the outer layer of the PCA. The control scheme is shown in Figure \ref{fig:ControlScheme}.

\subsubsection*{$\theta$ Tracking Law}
The dynamics of $\theta$ are, due to the inner layer of the PCA, simplified into a second-order integrator combined with sensor and actuator dynamics and their saturation. The sensor and actuator dynamics response time is on the order of milliseconds. Therefore, we do not include these dynamics in the control law and simplify the controlled dynamics into a second-order integrator system with saturation on control action. The control law comprises two hierarchical layers. The inner layer controls $\dot{\theta}$, while the outer layer controls the state $\theta$. Both control laws are simple proportional regulators with gains $G_{\dot{\theta}}$ and $G_\theta$, respectively. The controlled system has response dynamics
\begin{equation}
    \frac{\theta}{\theta_{\text{ref}}} = \frac{G_{\dot{\theta}}\cdot G_\theta}{s^2 + G_{\dot{\theta}}\cdot s + G_{\dot{\theta}}\cdot G_\theta}.
    \label{eq:CL_pitch}
\end{equation}

The $\dot{\theta}$ loop was parameterized to have a response time of 0.1 seconds and the $\theta$ loop 0.5 seconds. As the control law does not include integrators, no anti-windup methods needed to be implemented.

\subsubsection*{Lag Filter}
The outer layer of the PCA commands two variables: the pitch angle $\theta$ and the total thrust force $f_s$. The dynamics of $\theta$ are limited by the $\theta$ tracking law and can lag behind the $f_s$ commands in rapid reference changes. This behavior generates undesired behavior, typically overshoots in $\epsilon$. To prevent this behavior, we implement a lag filter on the variable $f_s$ with the same dynamical response as the $\theta$ tracking law.

\subsubsection*{$\epsilon$ Tracking Law}
The $\epsilon$ dynamics are modeled as a second-order integrator with actuation dynamics, corresponding to the $\theta$ tracking system and lag filter. We must emphasize that there is no direct translation between the higher- and lower-level linear dynamics. This is caused by the nonlinear dependency of $\ddot{\epsilon}$ and $\ddot{\lambda}$ on the variables $f_s$ and $\theta$. For this reason, we do not have knowledge of all system states. However, we can estimate the actuation-related states with the Kalman filter and then control the system with full state feedback. For details, refer to the provided git repository \cite{iuITwfwC9BcxK0Lg}.

\subsubsection*{$\lambda$ Tracking Law}
Identical to $\epsilon$ dynamics, the $\lambda$ dynamics are modeled as a second-order integrator with second-order actuation dynamics. We design a linear model predictive controller (MPC) to address large changes in $\lambda$ setpoints (e.g. multiple turns). Due to the PCA layer, which translates the nonlinear dynamics into linear ones, the designed linear MPC has consistent performance throughout the entire system envelope. We argue that the combination of linear MPC with a PCA layer is superior to the direct nonlinear MPC approach. Linear MPC has significantly smaller computational requirements, allowing for an increased prediction horizon and refresh rate with the same computing power. The design details are again available in the git repository \cite{iuITwfwC9BcxK0Lg}.

\section{Experimental Verification of PCA: Results}
\label{sec:Experiments}
The effectiveness and real-time applicability of the proposed Polynomial Control Allocation (PCA) method are validated through the following series of experiments on the Quanser 3DOF helicopter platform. The experiments were designed to evaluate the PCA method's ability to handle complex, nonlinear actuator dynamics and to demonstrate its superiority over traditional control allocation techniques in terms of precision, robustness, and computational efficiency.

The following sections detail the experimental procedures, present the results of step-change and multiple-turn tracking tasks, and analyze the computational demands observed during real-time operation.

\subsection{$\lambda$ and $\epsilon$ Step Change Tracking}
\begin{figure}[h]
    \centering
    \includegraphics[width=0.8\linewidth]{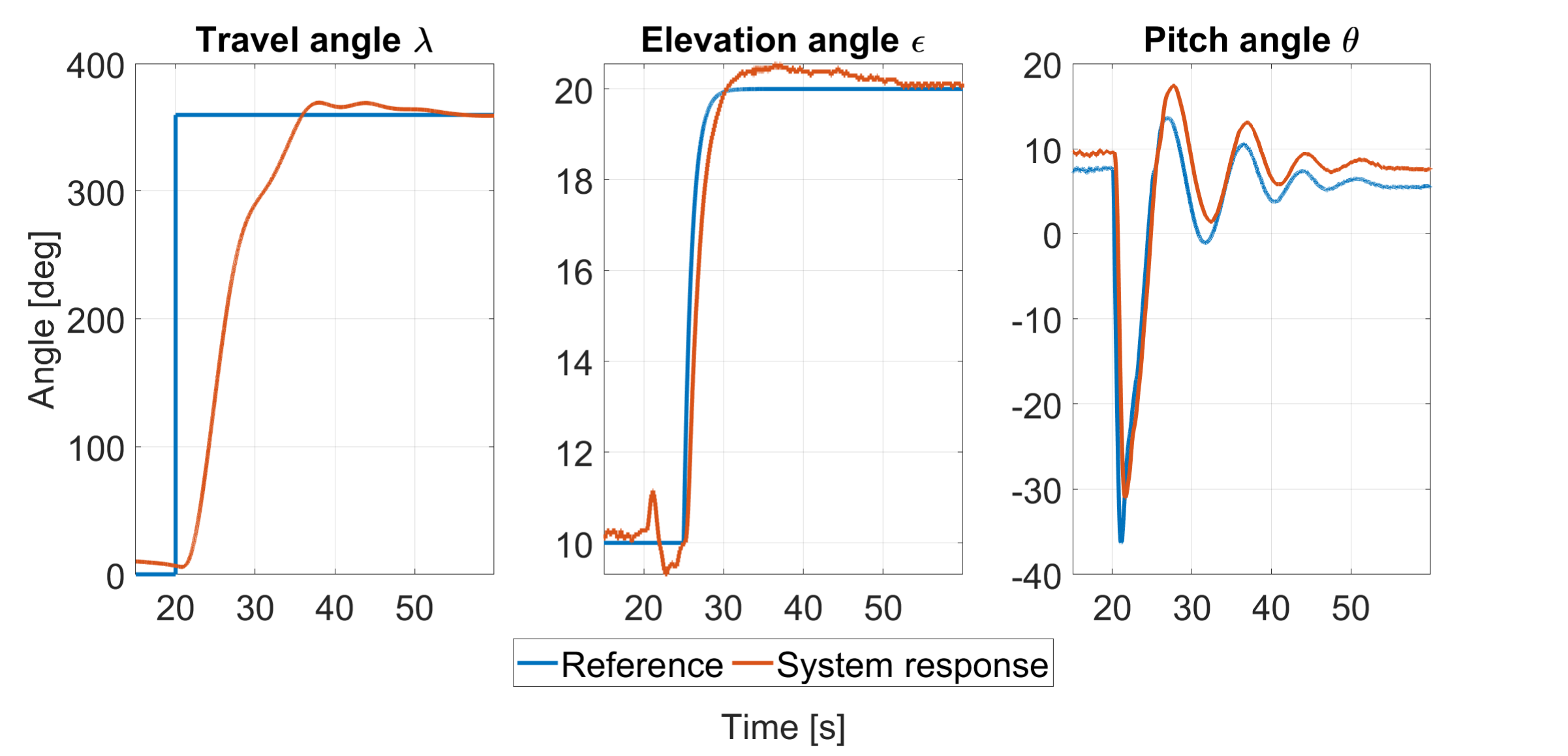}
    \caption{Step response to travel and elevation angles}
    \label{fig:experiment:Epsion_Lambda}
\end{figure}

In this experiment, step changes were applied independently to the reference commands of travel $\lambda$ and elevation $\epsilon$ of the controlled system. The primary objective was to demonstrate the effective decoupling of these two controlled variables using the proposed control allocation approach. The results are shown in Figure \ref{fig:experiment:Epsion_Lambda}.

As can be seen, the dynamics of $\lambda$ and $\epsilon$ are well decoupled. The 360-degree $\lambda$ reference step change causes less than a 1-degree error in $\epsilon$, which is quickly compensated by the $\epsilon$ controller. The step in the $\epsilon$ reference has no noticeable effect on the travel angle. The pitch angle $\theta$ is well tracked with the delay between the reference and the response of the system corresponding to the desired closed-loop dynamics (see equation \ref{eq:CL_pitch}.

\subsection{Multiple Full Turns in $\lambda$ Tracking}
\begin{figure}[h]
    \centering
    \includegraphics[width=0.8\linewidth]{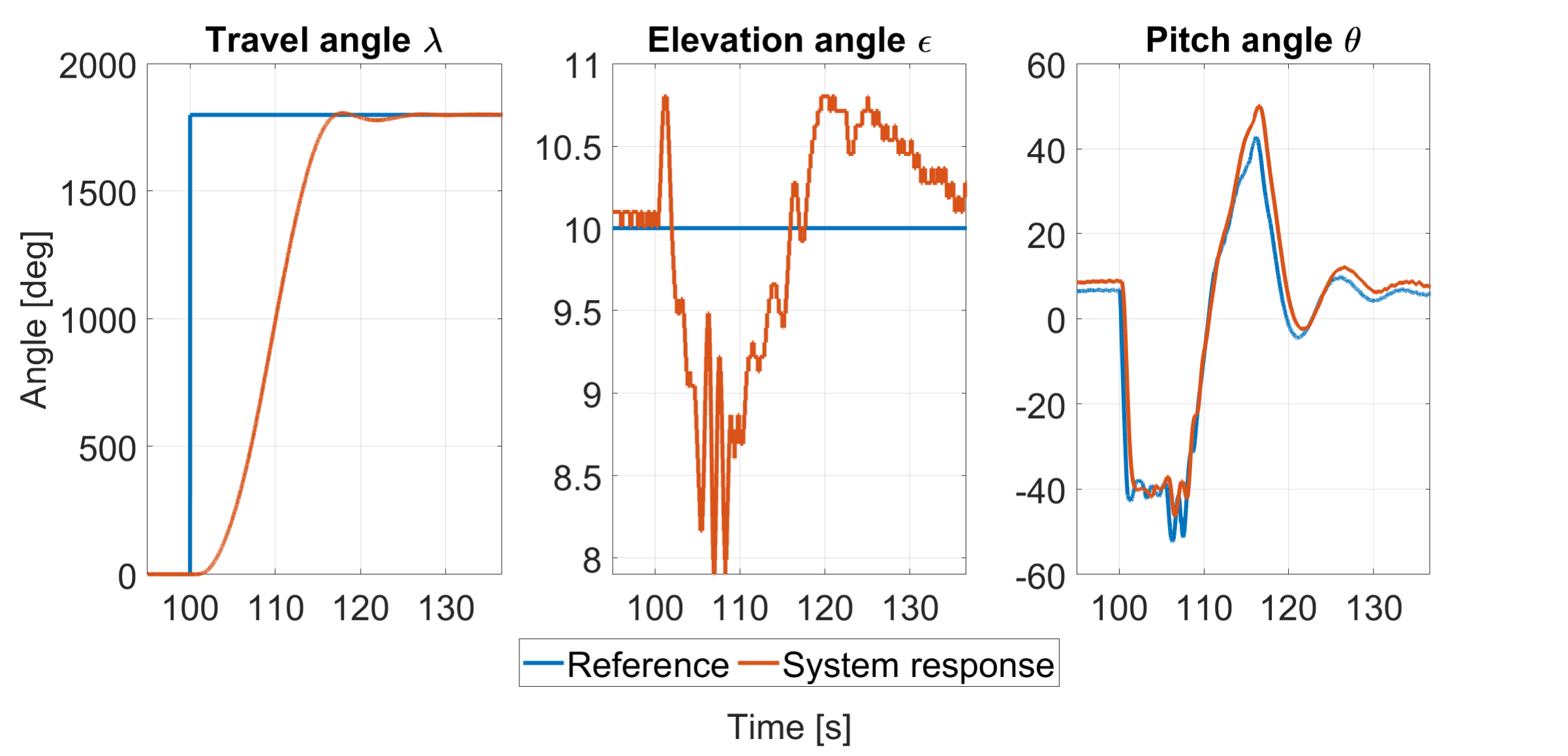}
    \caption{Multiple rotation of travel angle step response}
    \label{fig:experiment:Large_Lambda}
\end{figure}

This section presents experimental results from the tracking of multiple full rotations along the travel $\lambda$ axis. \hl{The focus is on evaluating system behavior within a nonlinear control region. The main objective is to demonstrate that a linear controller can operate outside the linearizable region of control authority, represented by large pitch angle $\theta$. Although the travel angle dynamics driven by pitch angle are highly nonlinear, the PCA method enables the controlled system to appear linear from the perspective of a higher-level controller, even at high pitch angles. Without this nonlinear control allocation approach, the designed linear control law (in our case linear MPC) would not be capable of operating in such conditions}.% While replacing it with a nonlinear MPC could allow operation in the nonlinear region without PCA, this would come at the cost of a significantly reduced prediction horizon and control frequency. A shorter prediction horizon would limit the controller’s ability to achieve the same maximum rotational speed, as there would be insufficient time for deceleration without overshoot. Therefore, the combination of PCA with linear MPC is argued to be comparable to, or potentially better than, nonlinear MPC in scenarios characterized by nonlinear control regions but linear autonomous system dynamics.} 
The resulting behavior is shown in Figure \ref{fig:experiment:Large_Lambda}.

As can be seen, the travel angle is tracked well with almost zero overshoot, which corresponds to the one-turn experiment from the previous section. The MPC reduces $\dot{\lambda}$ about five seconds before reaching the goal $\lambda$ through a large change in $\theta$. The experiment also demonstrated the decoupling of $\lambda$ and $\epsilon$ even for large $\theta$ values (up to 50 degrees), caused by large $\ddot{\lambda}$ commands at the beginning and end of the experiment.

\subsection{PCA Computational Loads}
\begin{figure}[h]
    \centering
    \includegraphics[width=0.7\linewidth]{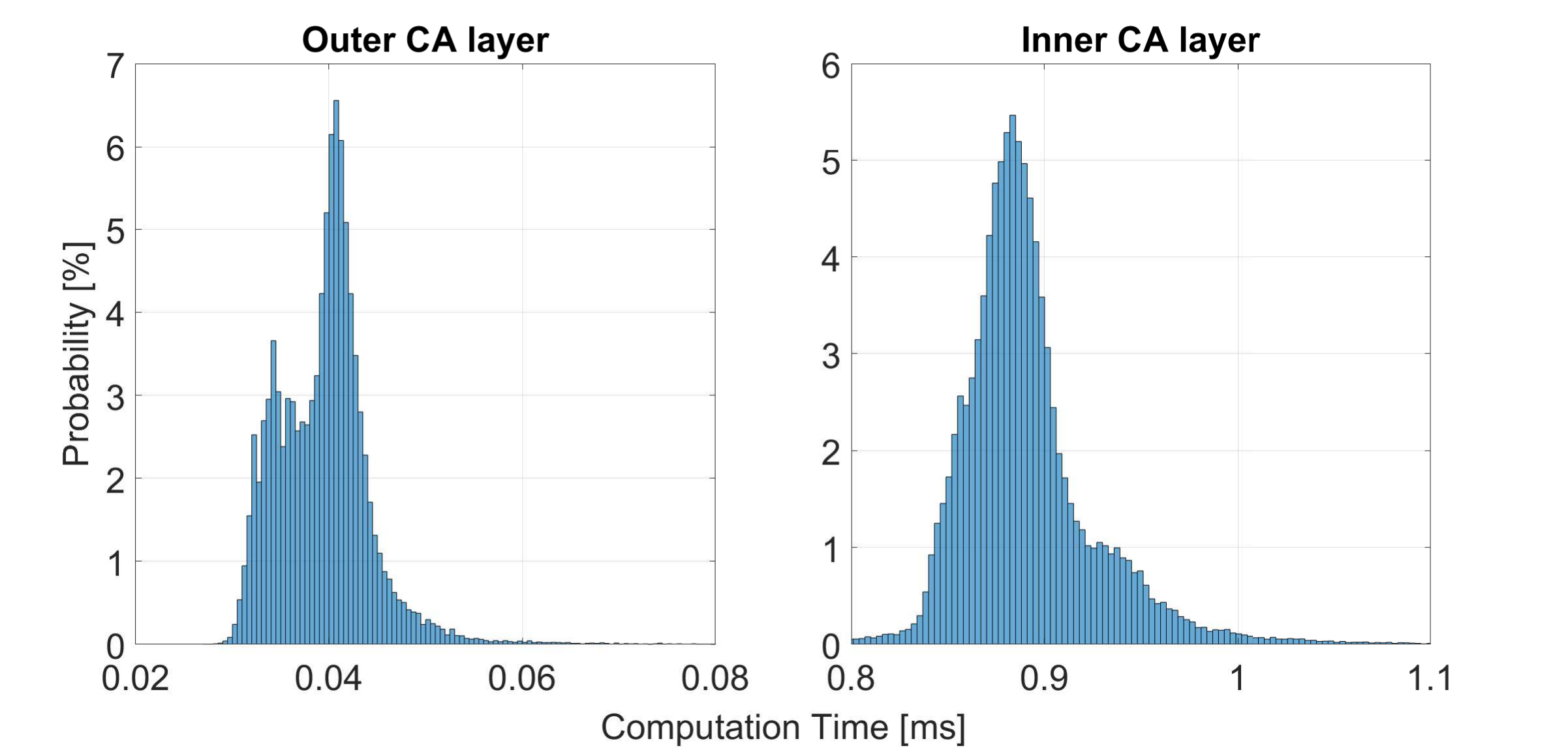}
    \caption{Computation time of PCA layers}
    \label{fig:experiment:PCA_comp_time}
\end{figure}

To verify the PCA computational complexity, we utilized simple MATLAB-based measurement of time needed to compute the PCA function. Each PCA level is examined independently to allow the reader to correlate computational complexity and problem difficulty. The histogram of computational time for the outer PCA layer is shown in Figure \ref{fig:experiment:PCA_comp_time}.

The typical computational time is about 0.05 ms, with the worst-case detected time of 0.25 ms. The worst case is possibly affected by the soft real-time design of the Quanser-MATLAB application. We argue that such a complex PCA layer can be implemented into systems with sampling periods as low as 0.5 ms, which is certainly sufficient for typical aircraft and multicopter control applications. 

%\begin{figure}[h]
%    \centering
%    \includegraphics[width=1\linewidth]{Low_CA_comp_hist.eps}
%    \caption{Computation time of inner PCA layer}
%    \label{fig:experiment:Low_PCA_comp_time}
%\end{figure}
The histogram of computational time for the inner PCA layer is shown in Figure \ref{fig:experiment:PCA_comp_time} on the right. Here, the nominal computational time is about 0.9 ms, with a worst case of about 1.4 ms. This means that in our application, the inner PCA layer utilizes about 45-70 percent of the available time. We assume that 1.4 ms is again caused by the Quanser real-time communication with the MATLAB application. The complexity of this layer should be suitable for 500 Hz or slower applications. These conclusions correspond well with the design expectations of 1 ms computational time, based on the number of solutions and polynomials in the matrices.

\section{Comparison to State of the Art}
\label{sec:Performance comparison to State of the art}
In the previous chapter, we described the performance of the PCA method on real hardware. With real-time capability and performance demonstrated, we focus in this chapter on comparison to the state-of-the-art methods described in Chapter \ref{sec:SOTA}. We utilize the upper CA layer of the Quanser control architecture described in Chapter \ref{sec:Quanser}. We compare the computational time and the relative CA error between the reference acceleration $a_r$ and the CA output $a_{\text{CA}}$, defined as
\begin{equation}
    \epsilon = \frac{a_{\text{CA}}-a_r}{a_r}.
\end{equation}

We compare the PCA method with two nonlinear CA methods from Chapter \ref{sec:SOTA}. The nonlinear programming (NPCA) method solves the programming problem described in the system of equations \ref{eq:OuterCA} using the sequential quadratic programming (SQP) algorithm. The Incremental nonlinear control allocation method (INCA) cannot solve this problem directly. Instead, it relies on local linearization which allows reformulating the problem as quadratic programming. To ensure that the linearized model is accurate in the set of possible solutions, we limit the deviation $d\theta$ and $dT$ from the previous value to be 5 degrees and 10\% of the full range, respectively, which means that it requires at least 36 steps to get from -180 degrees to 180 degrees. Furthermore, the previous values $\theta_0$ and $T_0$ are the measured/estimated values, not the output of the INCA method from the previous step. This ensures that the linearization is accurate even if the lower-level tracking law is slow. Both CA methods are implemented in the git repository \cite{iuITwfwC9BcxK0Lg}.

\subsection{Monte Carlo Evaluation}
To compare the PCA, NPCA, and INCA methods, we utilize $10^4$ Monte Carlo simulations with varying system state, described by $\epsilon$, $\dot{\epsilon}$, and $\dot{\lambda}$ variables, combined with varying acceleration requirements for $\ddot{\lambda}$ and $\ddot{\epsilon}$. The CA frequency is fifty times faster than the sample rate of Monte Carlo inputs to enable the INCA method to converge to the optimum. In these simulations, the INCA method utilizes the previous step as the values $\theta_0$ and $T_0$. Furthermore, the convergence time is directly proportional to the speed of the lower-level dynamics.

\begin{figure}[h]
    \centering
    \includegraphics[width=0.8\linewidth]{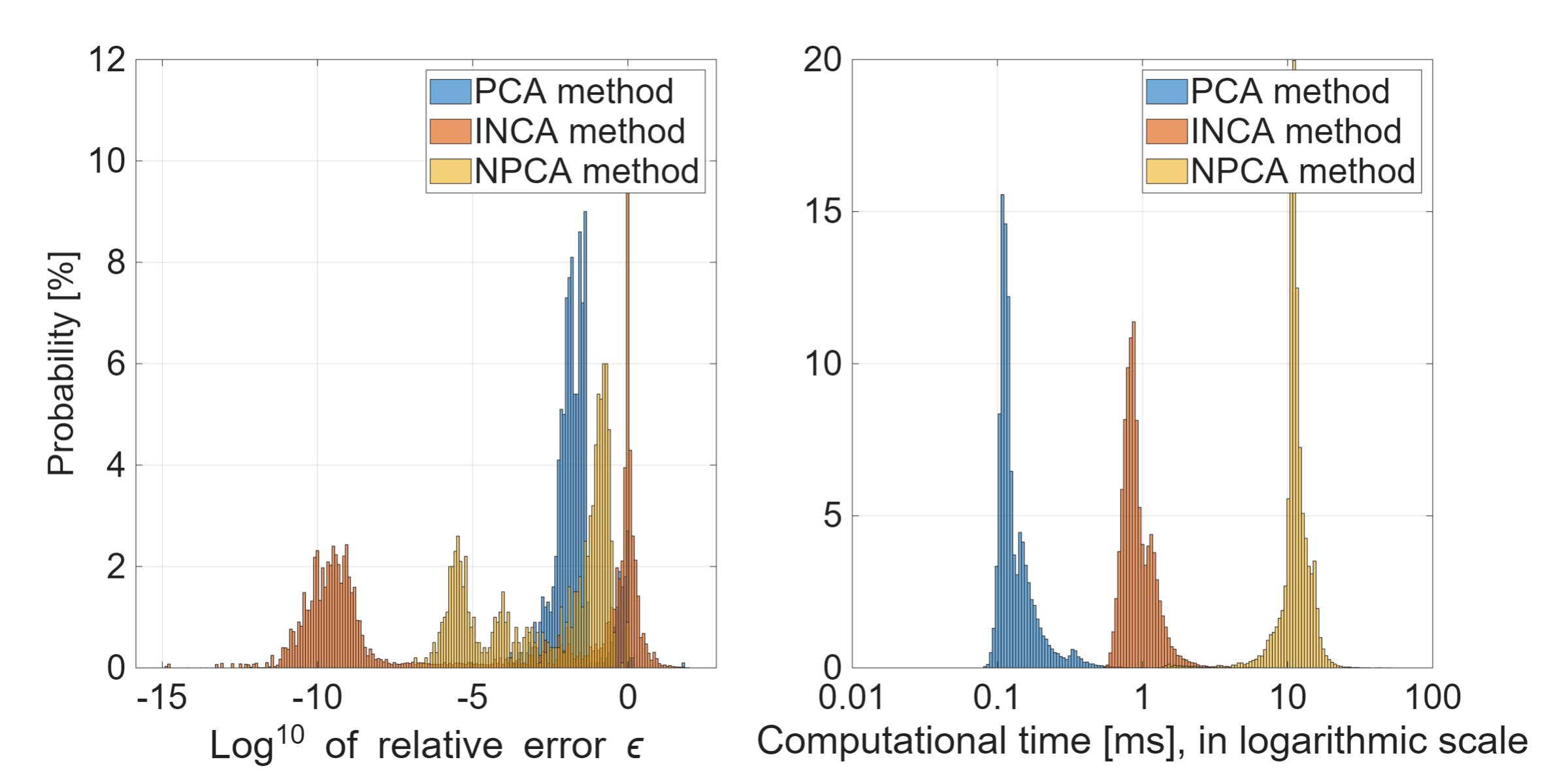}
    \caption{Comparison of nonlinear CA methods error rate (left) and computational time (right) using Monte Carlo simulations}
    \label{fig:comp_CA_error}
\end{figure}

The comparison of the resulting relative errors is shown in Figure \ref{fig:comp_CA_error} on the left. The PCA method is the most consistent, reliably achieving about 1\% relative error. The INCA method error rate varies greatly due to the iterative process described above. Eventually, the method converges to practically zero relative error; however, the time needed to achieve the optimum makes this method ill-suited for hierarchical control schemes. The NPCA method is more straightforward; however, the SQP method seems to struggle with the non-convexity of the problem. As a result, the NPCA method has a typical error rate closer to 10\%, likely representing local minima instead of global ones.

The comparison of computational time is shown in Figure \ref{fig:comp_CA_error} on the right. The PCA method is the fastest with about 0.5 ms time, followed by INCA method with about 1 ms and the NPCA with 10 ms computational time, respectively. The slow computation speed of INCA and NPCA methods is most likely connected to the termination requirements of their gradient descent algorithms. The computation time can be reduced; however, at the price of reduced CA accuracy.

%\begin{figure}[h]
%    \centering
%    \includegraphics[width=\linewidth]{Comp_CA_methods_time.eps}
%    \caption{Comparison of nonlinear CA methods computational time using Monte Carlo simulations}
%    \label{fig:comp_CA_time}
%\end{figure}

\section{Conclusion}
\label{sec:Conclusion}
This paper introduced the Polynomial Control Allocation framework, a novel approach that advances control allocation for nonlinear, over-actuated systems by formulating actuator dynamics as polynomial equations. Using Gröbner basis techniques with complex computations performed offline enables efficient real-time solutions for highly nonlinear, nonconvex problems previously intractable by traditional methods.

The PCA method broadens the class of actuator dynamics suitable for precise control allocation, overcoming the limitations of linearization, root-finding, and convex optimization techniques. The practical applications demonstrated include battery power control, drone force allocation, and optimal thrust control for variable pitch propellers. Experimental tests on the Quanser 3DOF helicopter demonstrated that fast, accurate, and robust control can be achieved thanks to the proposed PCA approach.

Performance analysis confirms PCA's suitability for real-time control applications through its combination of offline Gröbner basis computations and efficient online evaluations. Overall, PCA represents a significant advancement in control system design, enabling reliable real-time allocation for complex nonlinear systems.

\section*{Appendix}
\appendix
\section{Gröbner Basis Algorithms}
\label{appendix:GrobnerBasis}
Finding a solution to a system of polynomial equations is equivalent to finding a point in the set $V := \{x \in \mathbb{C}^n : p_1(x) = \cdots = p_m(x) = 0\}$ where $p_k$ are given elements of the ring $K[x]$. They are polynomials in $x$, our input and lifting variables, whose coefficients are elements of the base field $K$, rational functions of our parameters. To this geometric object $V$, called an algebraic variety, corresponds an algebraic object called an ideal and denoted $I$. It is the set of all polynomials that can be written as linear combinations $a_1\cdot p_1+\cdots+a_m\cdot p_m$ for some polynomial coefficients $a_k$. We say that polynomials $p_k$ are a representation of ideal $I$. Although $I$ has an infinite number of elements, a Hilbert theorem states that it always has a finite presentation, so that $I$ can be encoded by a computer. Note that every polynomial in $I$ vanishes at points of $V$, and we say that $V$ is the variety associated with $I$.

A key idea of algebraic geometry is to obtain useful information on $V$ from a suitable presentation of $I$. By taking finitely many linear combinations of polynomials $p_k$, we can generate another equivalent system of polynomials which generates the same ideal $I$ but with another presentation, and which is associated with the same variety $V$. In particular, when $V$ is a discrete set, i.e., a union of finitely many points, this presentation should allow computing the solutions easily. A useful presentation is a Gröbner basis. It can be computed with linear algebra and by performing a series of multivariate polynomial divisions. These divisions can be carried out provided one defines a suitable ordering on the set of monomials of variables $x$. To accelerate computations, these operations can be carried out in modular arithmetic, using appropriately large prime numbers.

If $V$ is discrete, then $I$ is a zero-dimensional ideal, and the quotient ring $K[x]/I$ is a finite-dimensional algebra. It is generated by finitely many monomials in $x$. The multiplication by each coordinate $x_k$ defines a linear map from $K[x]/I$ to itself, called a multiplication matrix denoted $M_k$, with entries in $K$. The values of $x_k$ on $V$ are the eigenvalues of $M_k$.

With MAPLE, these computations can be carried out using the Gröbner package. A Gröbner basis is computed with the `Basis` function. The monomials generating the quotient ring $K[x]/I$ are computed with the `NormalSet` function. The multiplication matrices are computed with the `MultiplicationMatrix` function. Beyond these general purpose commands, we refer the reader to \cite{Berthomieu2021718} and msolve.lip6.fr for state-of-the-art Gröbner basis software.

The multiplication matrices $M_k$ are computed once and for all. Their entries are explicit rational functions of the parameters. Then, for each value of the parameters, the multiplication matrices are evaluated numerically. A random linear combination $M(a) = a_1\cdot M_1 + \cdots + a_n\cdot M_n$ has distinct unit eigenvectors $v_1, \ldots, v_N$ such that $(v_l^T M_k\cdot v_l)_{k=1,\ldots,n} \in \mathbb{C}^n$ is a point of $V$, i.e., a solution $x$ of our polynomial equations, for $l = 1,\ldots,N$. From these $N$ complex solutions, we then sort out only those with zero imaginary parts and which are physically relevant.
\section*{Acknowledgment}
This research was supported by the Grant Agency of the Czech Technical University in Prague,
grant No. SGS25/145/OHK3/3T/13.

\end{document}